\documentclass{article}
\usepackage{scrextend}
\usepackage{geometry}
\usepackage[english]{babel}
\usepackage{amssymb}
\usepackage{tcolorbox}
\usepackage{amsmath}
\usepackage{amsthm}
\usepackage[normalem]{ulem}
\usepackage{float}
\useunder{\uline}{\ul}{}

\newtheorem{example}{Example}

\title{Introduction to Concentration Inequalities}
\date{}
\author{Kumar Abhishek \\\small IIIT Hyderabad, India.\\\small kumar.abhishek@research.iiit.ac.in 
\and Sneha Maheshwari\\\small IIT Roorkee, India. \\\small smaheshwari@ma.iitr.ac.in
\and Sujit Gujar\\\small IIIT Hyderabad, India.\\\small sujit.gujar@iiit.ac.in}
\usepackage[utf8]{inputenc}
\begin{document}
\maketitle 
\tableofcontents
\newpage
\begin{abstract}
In this report, we aim to exemplify concentration inequalities and provide easy to understand proofs for it. Our focus is on the inequalities which are helpful in the design and analysis of machine learning algorithms.  
\end{abstract}
\section{What are Concentration Inequalities?}
\emph{Concentration inequalities} furnish us bounds on how random variables deviate from a value (typically, expected value) or help us to understand how well they are concentrated. A random variable with high concentration is one that is close to its mean (or value) with high probability (more than a certain threshold). For example, the strong law of large numbers or weak law of large numbers say that under mild conditions, if we sum a large number of independant random variables,  with high probability, the sum is close to the expected value. These are elementary examples of the concentration we are talking about here. 
\par 
Concentration inequalities quantify the statements of random fluctuations of functions of random variables, typically by bounding the probability that such a function differs from its expected value (or from its median) by more than a certain amount.\\ \\
In the last decades, many researchers in a variety of areas were thriving to define concentration inequalities because of their importance in numerous applications.
\par
This report is organized as follows. In Section \ref{sec:motivation}, we provide examples where such concentration inequalities are useful. In Section \ref{sec:ineq}, we state and prove, (i) Markov's Inequality, (ii) Chebyshev's Inequality, (iii) Chernoff Bound, (iv) Hoeffding's Lemma, (v) Hoeffding's Inequality, (vi) Azuma's Inequality. In Section \ref{adineq}, we state advanced topics about concentration inequalities, (i) Bennett's Inequality, (ii) Bernstein's Inequality, (iii) Efron-Stein Inequality, (iv) McDiarmid's Inequality.   

\section{Motivation}
\label{sec:motivation}
\subsection{General example}
Let's start with a simple example,
\\ 
\begin{tcolorbox}
\textbf{Problem: Estimation of probability for a biased coin}
\\
Given a biased coin having an unknown probability '$p$' of occurring head, we need to estimate the value of $p$.
\end{tcolorbox}

\begin{itemize}
\item If we toss the coin once if it comes head then the probability of head will be $1$. But we are not at all confident for the probability being $1$.
\item If we toss the coin 100 times and head appears 65 times then we are a bit more confident for the probability being $0.65$.
\item Similarly, if we toss a coin lets say million times and the head is outcome $6,00,000$ times then we can say that '$p$' is $0.60$ with a very high confidence. 
\end{itemize}
Thus, to quantify the level of confidence with respect to the number of  trials, we can use \textit{concentration inequalities} to have better estimates of '$p$'. 
\subsection{Statistics}
In statistics we umpteen applications of concentration inequalities, let's see one of the example,
\begin{tcolorbox}
\textbf{Problem: Estimation of the population parameter.}
\\
In statistics, from an unknown population distribution, we want to infer information through sampling. (For example, one might want to know the population mean of age with probability of empirical mean to be actual mean, etc.)
\end{tcolorbox}
Following are the questions we need to address:
\begin{itemize}
\item How can we estimate the confidence interval (range of values) which would be a good estimate?
\item How can we determine the level of significance (confidence level) of that estimate?
\end{itemize}
We can get the answers to both questions through \textit{concentration inequalities}.
\subsection{Algorithms}
Zillions of analyses in algorithms (mainly in  theoretical computer science) uses concentration inequalities to give upper or lower bounds about the performance of algorithms with a certain probability. \\
For example, 
\begin{itemize}
\item \textbf{MAX cut problem}: We can solve this problem approximately and to analyze the probability that this algorithm gives a maximum cut we can use Reverse Markov inequality (converse of Markov inequality).
\end{itemize}
\subsection{Machine Learning}
In Machine learning, concentration inequalities are profoundly used in analyzing different aspects of learning algorithms. For example,
\begin{itemize}
\item Multi-Armed bandits problem: We use concentration inequalities to analyze  algorithms such as UCB algorithms, Thompson Sampling for their \emph{regret}, a measure on performance of a MAB algorithm. Here we need good estimates of rewards of each arm with high confidence. 
\end{itemize}
\subsection{Miscellaneous}
Among the areas of applications, without trying to be exhaustive, 
\begin{itemize}
\item Statistics
\item Learning Theory which includes  supervised learning, unsupervised learning, online learning, and reinforcement learning.
\item Discrete mathematics
\item Statistical mechanics
\item Information theory
\item High-Dimensional geometry
\end{itemize}
and the list goes on. In the next section, we will prove important concentration inequalities and illustrate with examples.  
\section{Inequalities}
We begin with the most elegant, yet powerful Markov inequality. Then, we go on explaining Chebyshev's inequality, Chernoff bound, Hoeffding's Lemma and inequality. At the end of this section, we state and prove Azuma's inequality. 
\label{sec:ineq}
\subsection{Markov's Inequality}
\begin{tcolorbox}
For a positive random variable $X \geq 0$ and $a > 0$, the probability that $X$ is no less than $a$ is less than or equal to the expectation of $X$ divided by $a$:
\begin{center}
$Pr[ X \geq a] \leq \frac{E(X)}{a}$
\end{center}
\end{tcolorbox}
\vspace{0.5cm}
\begin{proof}
\begin{equation}  \label{eq:1}
\begin{split}
E[X] & = \int_{0}^{\infty}xp(x)dx = \int_{0}^{a}xp(x)dx + \int_{a}^{\infty}xp(x)dx\\
& \geq \int_{a}^{\infty}xp(x)dx \geq a\int_{a}^{\infty}p(x)dx \\
& \geq aPr(X \geq a)
\end{split}
\end{equation}
By rearranging the terms, 

\begin{equation*}
Pr[ X \geq a] \leq \frac{E(X)}{a}    
\end{equation*}
\end{proof}
\begin{example}
Let $R$ be the weight distribution of a population with $E[R] = 100$. Calculate the probability that a random person weigh at least $200$ pounds.
\end{example}
\textbf{Solution:} As weight is always positive, we can apply Markov's inequality,\\
\begin{center}
$Pr[R \geq 200] \leq \frac{100}{200} \leq \frac{1}{2}$
\end{center}
\textbf{Corollary: Reverse Markov inequality}
\begin{tcolorbox}
Given maximum value 'U' of a random variable 'X',
\begin{center}
$Pr[X \leq a] \leq \frac{U - E[X]}{U - a}$
\end{center}
\end{tcolorbox}
Note: In the corollary there is no need for the random variable 'X' to be positive.
\begin{proof}
\begin{equation} \label{eq:2}
\begin{split}
Pr[X \leq a] & = Pr[U - X \geq U - a]\\
 & \leq \frac{E[U - X]}{U-a} \hspace{1cm} \text{(Applying Markov's inequality)} \\
 & \leq\frac{U - E[X]}{U - a} \hspace{1cm} 
\end{split}
\end{equation}
\end{proof}
\begin{example}
 Let 'X' be the random variable denoting the marks of random student. Maximum marks possible is $100$ (U) and expected marks $75$. What is the probability that a random student scores $50$ or less?\\
 \textbf{Solution:} We can directly apply reverse Markov inequality,
\begin{center}
$Pr[X \leq 50] \leq \frac{100 - 75}{100 - 50} \leq \frac{1}{2}$
\end{center}
\end{example}
\begin{example}
Suppose we use Markov's inequality to bound the probability of obtaining
more than 3n/4 heads in a sequence of n fair coin flips. Let 
\begin{equation*}
X_i=\begin{cases}1 & \text{if the $i^{th}$ coin flip is head}\\0 & \text{otherwise}\end{cases}
\end{equation*}
and let $X = \sum_{i=1}^n X_i$ denote the number of heads in the n coin flips. Since $E[X_i] =
Pr(X_i = 1) = 1/2$, it follows that $E[X] =\sum_{i=1}^n E[X_i] = n/2$. Applying Markov's
inequality, we obtain
\begin{equation*}
P(X\geq3n/4)\leq \frac{E[X]}{3n/4}=\frac{n/2}{3n/4}=\frac{2}{3}
\end{equation*}
\end{example}
\textbf{Features:}
\begin{itemize}
\item Upside - This needs almost no assumptions about the random variable.
\item Downside - It gives weaker bounds.
\end{itemize}
Markov's inequality is generally used where the random variable is too complicated to be analyzed by more powerful \footnote{Powerful inequalities are those whose confidence level are higher for small confidence interval }inequalities.
\subsection{Chebyshev's Inequality}
\begin{tcolorbox}
For a random variable X expectation and variance should be finite, then $\forall a > 0$,
\begin{center}
$Pr(|X - E[X]| \geq a) \leq \frac{Var[X]}{a^2}$
\end{center}
\end{tcolorbox}
\begin{proof}
\begin{equation} \label{eq:3}
\begin{split}
Pr(|X - E[X]| \geq a) & = Pr[(X - E[x])^2 \geq a^2]\\
 & \leq \frac{E[(X - E[X])^2]}{a^2} \hspace{0.2cm} \text{(Applying Markov's inequality)} \\
 & = \frac{Var[X]}{a^2} \hspace{1cm} 
\end{split}
\end{equation}
\end{proof}
\begin{example}
Let X be the IQ of random variable with $X \geq 0$, $E[X] = 100$ and $\sigma(X) = 15$. What is the probability of a random person having an IQ of atleast $250$?\\ \\
\textbf{Solution: }Let us first calculate using Markov's inequality,
\begin{center}
$Pr[X \geq 250] \leq \frac{100}{250} \leq 0.4$
\end{center}
Using Chebyshev's inequality we get, 
\begin{center}
$Pr[X - 100 \geq 150] \leq \frac{15^2}{150^2} \leq 0.01$
\end{center}
We can clearly see the difference on the bounds we got from the two concentration inequalities.\\
\end{example}
\begin{example}
Let us consider the coin-flipping example, and use Chebyshev's inequality to bound the probability of obtaining more than $3n/4$ heads in a sequence of n fair coin flips. Recall that $X_i = 1$ if the $i^{th}$ coin flip is heads and 0 otherwise, and $X = \sum_{i=1}^n X_i$ denotes the number of heads in the $n$ coin flips. To use Chebyshev's
inequality we need to compute the variance of $X$. Observe that, since $X_i$ is a bernoulli random variable,
\begin{equation*}
E[(X_i)^2] = E[X_i] = \frac{1}{2}
\end{equation*}
Thus,
\begin{equation*}
Var[X_i]= E[(X_i)^2]-(E[X_i])^2=\frac{1}{4}
\end{equation*}
Now, since $X = \sum_{i=1}^n X_i$ and the $X_i$ are independent
\begin{equation*}
Var[X]=Var[\sum_{i=1}^nX_i]=\sum_{i=1}^nVar[X_i]=\frac{n}{4}
\end{equation*}
Applying Chebyshev's inequality yields
\begin{equation*}
\begin{split}
P[X\geq3n/4]&\leq P[|X-E[X]|\geq n/4]\\
&\leq \frac{Var[X]}{(n/4)^2}\\
&=\frac{(n/4)}{(n/4)^2}\\
&=\frac{4}{n}
\end{split}
\end{equation*}
In fact, we can do slightly better. Chebyshev's inequality yields that $4/n$ is actually
a bound on the probability that $X$ is either smaller than $n/4$ or larger than $3n/4$, so
by symmetry the probability that $X$ is greater than $3n/4$ is actually $2/n$. Chebyshev's
inequality gives a significantly better bound than Markov's inequality for large n.
\end{example}
\textbf{Usage: }Chebyshev's inequality has great utility because it can be applied to any probability distribution in which the mean and variance are defined.
\subsection{Chernoff Bound}
\begin{tcolorbox}
The generic Chernoff bound for a random variable X states,
\begin{center}
$Pr(X \geq a) = Pr(e^{tX} \geq e^{ta}$) \hspace{0.2cm} $\forall t > 0$
\end{center}
\end{tcolorbox}
As $e^{tX} \geq 0$ and is monotonically increasing function, we can use Markov's inequality,
\begin{center}
$Pr(X \geq a) \leq \frac{E[e^{tX}]}{e^{ta}}$
\end{center}
When $X = X_1 + X_2 ....+ X_n$ for any $t > 0$,
\begin{center}
$Pr(X \geq a) \leq e^{-ta} \hspace{0.25cm}E\big[\prod_i e^{tX_i}\big]$
\end{center}
For better tighter bounds we can optimize over '$t$'.\\
\\
\textbf{Derivation of Chernoff bound for Bernoulli random variable}\\
Let $X_1, X_2, ......,X_n$ be independent rv(random variable), whose sum is $X$.\\
Let '$p$' be the probability of $X_i = 1$.
\begin{equation} \label{eq:4}
\begin{split}
E[e^{tX_i}] & = pe^t +(1-p)\\
& = 1 + p(e^t - 1)\\
& \leq e^{p(e^t - 1)}\hspace{1cm} (1 + x \leq e^x)
\end{split}
\end{equation}
\begin{equation} \label{eq:5}
\begin{split}
Pr(X \geq a) & \leq \frac{E[e^{tX}]}{e^{at}}\\
& \leq e^{-at}E[e^{\sum_i tX_i}]\\
& \leq e^{-at}E[e^{tX_1}]\hspace{0.01cm}[e^{tX_2}]....E[e^{tX_n}]\\
& \hspace{0.3cm}\text{(As given independent rv's)}\\
& \leq e^{-at}e^{\sum_ip(e^t - 1)} \hspace{0.5cm}\text{(From Eq. (\ref{eq:4}))}
\end{split}
\end{equation}
Now, substitute the following for $\delta > 0$ in Eq. (\ref{eq:5}),
\begin{equation} \label{eq:6}
\begin{split}
a & = (1 + \delta)np\\
& = (1 + \delta)E[X]\\
t & = \ln(1 + \delta)
\end{split}
\end{equation}
We will get,
\begin{equation} \label{eq:7}
\begin{split}
Pr(X \geq (1 + \delta)np) & \leq \frac{e^{np(1 + \delta - 1)}}{(1 + \delta)^{(1 + \delta)np}}\\
& \leq [\frac{e^{\delta}}{(1 + \delta)^{1 + \delta}}]^{np}
\end{split}
\end{equation}
Similarly, we can derive for different random variables.
\begin{example}
\label{example:}
$1$ million people are playing pick $4$ ($0000 - 9999$), i.e., there is a fixed $4$ digit number and all people have to guess the number to be the winner. Calculate the probability of atleast $200$ winner's.
\begin{center}
$Pr[win] = \frac{1}{10000}$\\
$E[\text{Number of winners}] = 100$.\\
\end{center}
\begin{equation*}
\begin{split}
Pr[X \geq 200] & = Pr[X \geq (1+\delta)100] \hspace{0.25cm}(\text{where } \delta = 1) \\
& \leq [\frac{e}{2^2}]^{100} \\
& \leq (0.67)^{100} = 4.05 * e^{-18}
\end{split}
\end{equation*}
We got a very small probability, hence we have a tight bound.
\end{example}
\begin{example}
Let X be the number of heads in a sequence of n independent fair coin flips.To compare the power of this bound to Chebyshev's bound. consider the probability
of having no more than n/4 heads or no fewer than 3n/4 heads in a sequence of n
independent fair coin flips. In the previous theorem, we used Chebyshev's inequality to
show that
\begin{equation*}
P\bigg(\bigg|X-\frac{n}{2}\bigg|\geq\frac{n}{4}\bigg)\leq\frac{4}{n}
\end{equation*}
Using the Chernoff bound in this case, we find that
\begin{equation*}
P\bigg(\bigg|X-\frac{n}{2}\bigg|\geq\frac{n}{4}\bigg)\leq 2exp\bigg\{-\frac{1}{3}\frac{n}{2}\frac{1}{4}\bigg\} = 2e^{-\frac{n}{24}}
\end{equation*}
Observe that Chernoff bound gives a bound that is exponentially smaller than the bound obtained using Chebyshev's inequality.
\\
\end{example}
\textbf{Applications: }
\begin{itemize}
\item Chernoff bound is used to bound the tails of the distribution for a sum of independent random variables.
\item The Chernoff bound is by far the most useful tool in randomized algorithms.
\item Application in Networking : Chernoff bound is also used to obtain tight bounds for permutation routing problems which reduce network congestion while routing packets in sparse networks.
\end{itemize}
\vspace{0.5cm}
Summarizing the above three inequalities,
\begin{itemize}
\item Markov's Inequality : This inequality suffices when constant probability bound is sufficient for the task.
\item Chebyshev's Inequality : This inequality is the appropriate one when one have a good handle on the variance of the random variable.
\item Chernoff bound : This inequality gives sharp concentration bounds for random variables that are sums of independent and bounded random variables (most commonly,
sums of independent indicator random variables).
\end{itemize}
\subsection{Hoeffding's Lemma}
Hoeffding's lemma is an inequality that bounds the moment-generating function of any bounded random variable.\par
Note that Markov's inequality bounded first moment of random variable and Chebyshev's bounded second moment of random variable.
\vspace{0.5cm}
\begin{tcolorbox}
Let $X$ be any real valued random variable with $E[X] = \mu$, such that $a\leq X \leq b$ almost surely (that is with probability $= 1$). Then $\forall \lambda \in R$,
\begin{center}
$E[e^{\lambda X}] \leq e^{\lambda \mu} e^{\frac{(\lambda)^2 (b - a)^2}{8}}$
\end{center}
\end{tcolorbox}
\begin{proof}
As exponential function in convex we will use convexity property, we can write $X$ as convex combination of $a$ and $b$.
\begin{equation}\label{eq:8}
\begin{split}
X & = tb + (1 - t)a \hspace{0.5cm} \text{where } t \in [0,1]\\
t & = \frac{X -a}{b - a}
\end{split}
\end{equation}
\vspace{0.5cm}
\begin{equation} \label{eq:9}
\begin{split}
e^{\lambda X} & = e^{\lambda(tb + (1 - t)a)} \hspace{3cm} \\
& \leq te^{\lambda b} + (1 - t)e^{\lambda a}
\end{split}
\end{equation}
Taking expectation and substituting Eq. (\ref{eq:8}) in Eq. (\ref{eq:9}),
\begin{equation}\label{eq:10}
\begin{split}
E[e^{\lambda X}] & \leq e^{\lambda b}\hspace{0.01cm}E\bigg[\frac{X - a}{b - a}\bigg] +e^{\lambda a}\hspace{0.01cm}E\bigg[\frac{b - X}{b - a}\bigg]\\
& \leq  e^{\lambda b}\hspace{0.01cm}\bigg(\frac{\mu - a}{b - a}\bigg) +e^{\lambda a}\hspace{0.01cm}\bigg(\frac{b - \mu}{b - a}\bigg)\\
& \hspace{0.5cm}\bigg(\text{Now substituting $ \gamma = \frac{b - \mu}{b - a}$ }\bigg)\\
& \leq e^{\lambda b}(1 - \gamma) + e^{\lambda a}\gamma
\end{split}
\end{equation}
Let $u = (b - a)\lambda$. Consider the following function:
\begin{equation}\label{eq:11}
\begin{split}
\phi(u) & = \log(\gamma e^{\lambda a} + (1 - \gamma)e^{\lambda b})\\
& = \lambda a + \log((1 - \gamma)e^u + \gamma)\\
& = (\gamma - 1)u + \lambda \mu + \log((1 - \gamma)e^u + \gamma)
\end{split}
\end{equation}
As, $E[e^{\lambda X}] \leq e^{\phi (u)}$. To find the least upper bound, we need to minimize $\phi(u)$.
\par
$\phi(u)$ is twice differentiable and hence using Taylor's theorem for any $u$ there exists $\xi \in [0, u]$ such that,
\begin{equation*}
\phi(u) = \phi(0) + u{\phi}^{'}(0) + \frac{u^2}{2} \phi^{''}(\xi)
\end{equation*}
Using Eq. (\ref{eq:11}), we can see that $\phi(0) = \lambda \mu$. Also, 
\begin{equation}\label{eq:12}
\begin{split}
{\phi}^{'}(0) & = (1 - \gamma) + \frac{(1 - \gamma)e^u}{\gamma + (1 - (\gamma)e^u}\\
\phi^{''}(u) & = \frac{(1 - \gamma)e^u}{\gamma + (1 - (\gamma)e^u}[ 1 - \frac{(1 - \gamma)e^u}{\gamma + (1 - (\gamma)e^u}]
\end{split}
\end{equation}
Thus, $\phi^{'}(0) = 0$, $\phi^{''}(u) =  p(1 - p)$, where $p = \frac{(1 - \gamma)e^u}{\gamma + (1 - (\gamma)e^u}$. Thus, $\phi^{''}(u) \leq \frac{1}{4}$. Hence,
\begin{center}
$\phi(u) \leq \lambda \mu + \frac{1}{8}u^2 = \lambda \mu + \frac{\alpha^{2}}{8}(b - a)^2$
\end{center}
Thus,
\begin{equation}\label{eq:13}
\begin{split}
E[e^{\lambda X}] & \leq e^{\phi(u)}\\
& \leq e^{\lambda \mu}e^{\frac{\lambda^{2}(b - a)^2}{8}}
\end{split}
\end{equation}
\end{proof}
\subsection{Hoeffding's Inequality}
 Hoeffding's inequality provides an upper bound on the probability that the sum of independent random variables deviates from its expected value by more than a certain amount.
 \begin{tcolorbox}
 Let $X_1, X_2,.....,  X_n$ be $n$ independent random variables, and $S_n = X_1 + X_2 +....+ X_n$, where $\forall i, X_i \in [a_i ,b_i]$, then according to Hoeffding's inequality,
 \begin{center}
 $Pr[S_n - E[S_n] \geq t] \leq e^{\frac{-2t^2n^2}{\sum_i (b_i - a_i)^2}}$
 \end{center}
 \end{tcolorbox}
 \begin{proof}
 \begin{equation}\label{eq:14}
 \begin{split}
 Pr[S_n - E[S_n] \geq t] & = Pr[e^{s(S_n - E[S_n])} \geq e^{st}] \hspace{0.5cm} (\text{For }\forall s > 0) \\
 & \leq \frac{E[e^{s[S_n - E[S_n]]}]}{e^{st}} \hspace{0.5cm} (\text{Applying Markov's inequality})
 \end{split}
 \end{equation}
 \begin{equation*}
 \begin{split}
 E[e^{s[S_n - E[S_n]]}] & = E[e^{s\sum_i^n X_i - E[X_i]}]\\
 & = E\bigg[\prod_i^n e^{s(X_i - E[X_i]}\bigg]\\
 & \hspace{0.5cm}(\text{Substituting } Y_i = X_i - E[X_i])\\
 \end{split}
 \end{equation*}
 \begin{equation} \label{eq:15}
 \begin{split}
 E[e^{s[S_n - E[S_n]]}] & = E\bigg[{\prod_i^n e^{s Y_i}}\bigg]\\
 & \leq \prod_i^n[e^{sE[Y_i]}e^{\frac{s^2(b_i - a_i)^2}{8}}] \hspace{0.5cm}(\text{Applying Hoeffding's Lemma})\\
 & \leq \prod_i^n e^{\frac{s^2(b_i - a_i)^2}{8}} \hspace{0.5cm}({E[Y_i] = 0})
 \end{split}
 \end{equation}
 By substituting Eq. (\ref{eq:15}) in Eq. (\ref{eq:14}), we get
 \begin{equation}\label{eq:16}
  Pr[S_n - E[S_n] \geq t] \leq e^{-st + \frac{s^2 \sum_i^n (b_i - a_i)^2}{8}}
 \end{equation}
 To get the best possible upper bound, we find the minimum of the right hand side of the last inequality as a function of $s$. Define
 \begin{equation*}
 g(s) = -st + \frac{s^2 \sum_i^n (b_i - a_i)^2}{8}
 \end{equation*}
 Note that $g$ is a quadratic equation and achieves its minimum at
 \begin{equation*}
 s = \frac{4t}{\sum_i^n (b_i - a_i)^2}
 \end{equation*}
 Thus we get 
 \begin{equation}\label{eq:17}
  Pr[S_n - E[S_n] \geq t] \leq e^{\frac{-2t^2}{\sum_i^n (b_i - a_i)^2}}
 \end{equation}
 \end{proof}
 \begin{flushleft}
 \textbf{Usage}
 \end{flushleft}

 One of the main application of Hoeffding's inequality is to analyse the number of required samples needed to obtain a confidence interval by solving the inequality,
 \begin{equation*}
 Pr[\bar{X} - E[\bar{X}] \geq t] \leq e^{-2nt^2}
 \end{equation*}
 Symmetrically, the inequality is also valid for another side of the difference:
 \begin{equation*}
 Pr[-\bar{X} + E[\bar{X}] \geq t] \leq e^{-2nt^2}
 \end{equation*}
 By adding them both up, we can obtain two-sided variant of this inequality:
 \begin{equation*}
 Pr[|\bar{X} - E[\bar{X}|] \geq t] \leq 2e^{-2nt^2}
 \end{equation*}
 This probability can be interpreted as the level of significance $ \alpha $(probability of making an error) for a confidence interval around $E[\bar{X}]$ of size $2t$:
 \begin{equation*}
 \alpha = P(\bar{X} \notin [E[\bar{X}] - t ,E[\bar{X}] + t]) \leq 2e^{-2nt^2}
 \end{equation*}
 Solving for the number of required samples $n$ gives us,
 \begin{equation*}
 n \geq \frac{log(2/\alpha)}{2t^2}
 \end{equation*}
 Therefore, we require at least $\frac{log(2/\alpha)}{2t^2}$ samples to acquire $(1 - \alpha)$ confidence interval $E[\bar{X}]\pm t$.
 \subsection{Azuma's Inequality}
 The Azuma–Hoeffding inequality gives a concentration result for the values of martingales that have bounded differences. That is here random variables are not independent.
 \begin{tcolorbox}
 Let $Z_0,...,Z_n$ be a martingale sequence with respect to the filter $\mathcal{F}_0 \subseteq \mathcal{F}_1 \subseteq .... \subseteq \mathcal{F}_n$ such that for $Y_i = Z_i - Z_{i-1}$, we have that for all $i \in [n]$, $|Y_i| = |Z_i - Z_{i-1}| \leq c_i$. Then 
 \begin{center}
 $Pr[Z_N - Z_0 \geq t] \leq exp\bigg(\frac{-t^2}{2\sum_{i = 1}^n c_i^2}\bigg)$ and $Pr[Z_0 - Z_n \geq t] \leq exp\bigg(\frac{-t^2}{2\sum_{i = 1}^n c_i^2}\bigg)$
 \end{center}
 \end{tcolorbox}
\begin{proof} 
We first prove one side of inequality. For any $\lambda > 0$, using Chernoff bound and Markov's inequality\\
 \begin{center}
 $Pr[Z_n - Z_0 \geq t] = Pr[e^{\lambda(Z_n - Z_0)} \geq e^{\lambda t}] \leq e^{-\lambda t} \mathbb{E}[e^{\lambda(Z_n - Z_0)}]$
 \end{center}
 Now conditioning on $\mathcal{F}_{n-1}$, we get
  \begin{equation*}
  \begin{split}
   \mathbb{E}[e^{\lambda(Z_n - Z_0)}] & = \mathbb{E}[e^{\lambda(Y_n + Z_{n-1} - Z_0)}] \\
   & = \mathbb{E}[\mathbb{E}[e^{\lambda(Y_n + Z_{n-1} - Z_0)} | \mathcal{F}_{n-1}]] \\
   & = \mathbb{E}[e^{\lambda(Z_{n-1} - Z_0)} \mathbb{E}[e^{\lambda Y_n} | \mathcal{F}_{n-1}]]
  \end{split}
 \end{equation*}
 Using the fact that $Z_{n-1}$ and $Z_0$ are both measurable in the $\sigma-$algebra $\mathcal{F}_{n-1}$. We not bound the expectation $\mathbb{E}[e^{\lambda Y_n} | \mathcal{F}_{n-1}]$ using convexity of the function $e^x$. Let $\alpha \in [-1, 1]$ and $M \in \mathbb{R}$ be any real number. Then, 
 \begin{equation*}
 \alpha M = \bigg(\frac{1 + \alpha}{2}\bigg)M - \bigg(\frac{1 - \alpha}{2}\bigg)M
 \end{equation*}
 Now using the convexity of the function $e^x$, 
 \begin{equation*}
 e^{\alpha M} \leq \bigg(\frac{1 + \alpha}{2}\bigg)e^M + \bigg(\frac{1 - \alpha}{2}\bigg)e^{-M}
 \end{equation*}
 Now taking $\alpha = Y_n/c_n$ and $M = \lambda c_n$, we get 
 \begin{equation*}
 e^{\lambda Y_n} \leq \Bigg(\frac{1 + (Y_n/c_n)}{2}\Bigg) e^{\lambda c_n} + \Bigg(\frac{1 - (Y_n/c_n)}{2}\Bigg) e^{-\lambda c_n}
 \end{equation*}
 Using $\mathbb{E}[Y_n | \mathcal{F}_{n-1}] = 0$, we get
 \begin{equation*}
 \begin{split}
 \mathbb{E}[e^{\lambda Y_n }| \mathcal{F}_{n-1}] & \leq \mathbb{E} \Bigg[\Bigg(\frac{1 + (Y_n/c_n)}{2}\Bigg) e^{\lambda c_n} + \Bigg(\frac{1 - (Y_n/c_n)}{2}\Bigg) e^{-\lambda c_n} \Bigg| \mathcal{F}_{n-1}\Bigg] \\
 & = \frac{e^{\lambda c_n }+ e^{- \lambda c_n}}{2} \leq e^{\frac{(\lambda c_n)^2}{2}}
 \end{split}
 \end{equation*}
 where the last step uses the fact $(e^x + e^{-x})/2 \leq e^{\frac{x^2}{2}}$ which uses taylor expansion to verify. 
 \begin{equation*}
 Pr[Z_n - Z_0 \geq t] \leq e^{-\lambda t} e^{\lambda^2 c_n^2 / 2} \mathbb{E}[e^{\lambda(Z_{n-1} - Z_0)}]
 \end{equation*}
 Continuing by same process, we can deduce
 \begin{equation*}
 Pr[Z_n - Z_0 \geq t] \leq exp\bigg(-\lambda t + ({\lambda}^2 / 2)\sum_{i =1}^n c_i^2 \bigg)
 \end{equation*}
 Since above equation holds for any $\lambda > 0$, we can optimize over $\lambda$ to minimize the above bound. On calculating the above expression is minimized for $\lambda = \frac{t}{\sum_{i =1}^n c_i^2} $, which gives
 \begin{equation*}
 Pr[Z_n - Z_0 \geq t] \leq exp\bigg( -\frac{t^2}{2\sum_{i =1}^n c_i^2}\bigg)
 \end{equation*}
 \end{proof}
\begin{flushleft}
 Similarly it can be proven for $Pr[Z_0 - Z_n \geq t]$.
\end{flushleft}
 \begin{example}
 Some times, we have to find the the interesting patterns, example examining DNA structure.

Let $X=(X_1,...,X_n)$ be independent characters chosen from alphabet $A$ where $a=|A|$. Let $B=(b_1,...,b_k)$ be fixed string of $k$ characters from $A$. Let $F$ be the number of occurrence of the fixed string $B$ in the random string $X$.\\
Let,
\begin{equation*}
Z_0=E[F]
\end{equation*}
and for $1\leq i \leq n$ let 
\begin{equation*}
Z_i=E[F|X_1,...,X_i]
\end{equation*}
The sequence $ Z_0,... , Z_n$ is a Doob martingale, and
\begin{equation*}
Z_n=F
\end{equation*}
Since each character in the string X can participate in no more than k possible
matches, for any $0 \leq i \leq n$ we have
\begin{equation*}
|Z_{i+1}-Z_i| \leq k
\end{equation*}
\\
In other words. the value of $X_{i+1}$ can affect the value of F by at most k in either direction,
since $X_{i+1}$ participates in no more than k possible matches. Hence the difference is
\begin{equation*}
E[F|X_1,...,X_{i+1}]-E[F|X_1,...,X_i]|=|Z_{i+1}-Z_{i}|
\end{equation*}
must be at most k, Applying Azuma-Hoeffding Inequality yields
\begin{equation*}
P[|F-E[F]|\geq\epsilon]\leq 2e^{\frac{-\epsilon^2}{2nk^2}}
\end{equation*}
 \end{example}
 
 \section{Advanced Inequalities}
 \label{adineq}
 In this section we now study advanced inequalities, namely:  Bennett's Inequality, Bernstein's Inequality, Efron-Stein Inequality, McDiarmid's Inequality.
 \subsection{Bennett's Inequality}
 \begin{tcolorbox}
 Let $X_1,...,X_n$ be independent real-valued
random variables with zero mean, and $|X_i| \leq 1$ with probability one. Then for any $t > 0$
\begin{center}
$\mathbb{P}\bigg[\sum_{i=1}^n X_i > t\bigg] \leq exp \bigg(-n\sigma^2 h \bigg(\frac{t}{n\sigma^2}\bigg)\bigg)$
\end{center}
 \end{tcolorbox}
 where,
 \begin{equation*}
 \sigma^2 = \frac{1}{n}\sum_{i =1}^n Var\{X_i\}
 \end{equation*}
 \begin{equation*}
 h(u) = (1 + u)\log(1 + u) - u \ for\  u \geq 0
 \end{equation*}
 \begin{proof}
 
Given that mean of rv's are zero , that is 
\begin{equation} \label{eq:18}
E[X_i]=0
\end{equation}
Let 
\begin{equation}\label{eq:19}
F_i =\sum_{r=2}^\infty\frac{s^{r-2}E(X_i^r)}{r!\sigma_i^2}
\end{equation}
where \space $ \sigma_i^2 = E(X_i^2)-E(X_i)^2 = Var\{X_i\}$
\\
\\
now , $ e^x = 1 + x + \sum_{r=2}^\infty\frac{x^r}{r!} $ therefore, 
\begin{equation}\label{eq:20}
\begin{split}
E(e^{sX_i})&=1+sE(X_i)+\sum_{r=2}^\infty\frac{s^rE(X_i^r)}{r!}\\
E(e^{sX_i}) &= 1+ s^2\sigma_i^2F_i  \hspace{1cm} \text{(Using Eq. (\ref{eq:18}) and Eq. (\ref{eq:19}))}\\
&\leq e^{ s^2\sigma_i^2F_i}
\end{split}
\end{equation}
Consider the term $E(X_i^r)$. Since expectation of a function is just the Lebesgue
integral of the function with respect to probability measure, we have \\
$E(X_i)=\int_{P} X_i^{r-1}X_i$.
\ Using Cauchy Schwarz inequality we get,
\begin{equation*}
\begin{split}
E(X_i^r)&=\int_{P}{X_i^{r-1}X_i}\\
& \leq\bigg(\int_{P}\mid X_i^{r-1}\mid^2\bigg)^{1/2} \bigg(\int_{P}\mid X_i\mid^2)^{1/2}\bigg)\\
\Rightarrow E(X_i^r)& \leq \sigma_i\bigg(\int_{P}\mid X_i^{r-1}\mid^2\bigg)^{1/2}
\end{split}
\end{equation*}
Proceeding to use the Cauchy Schwarz inequality recursively $k$ times we get
\begin{equation*}
\begin{split}
E(X_i^r) & \leq\sigma_i^{1+\frac{1}{2}+\frac{1}{2^2}+...+\frac{1}{2^{k-1}}}\bigg(\int_{P}\mid X_i^{(2^kr-2^{k-1}-1)}\mid\bigg)^{1/2^k}\\
& = \sigma_i^{2(1-\frac{1}{2^k})}\bigg(\int_{P}\mid X_i^{(2^kr-2^{k-1}-1\bigg)}\mid)^{1/2^k}\\
\end{split}
\end{equation*}
Now we know that $|X_i| \leq 1$. Therefore,
\begin{equation*}
\bigg(\int_{P}\mid X_i^{(2^kr-2^{k-1}-1)}\mid\bigg)^{1/2^k}\leq 1
\end{equation*}
Hence, we get 
\begin{equation*}
E(X_i^r)\leq \sigma_i^{2(1-\frac{1}{2^k})}
\end{equation*}
Taking limit $k \rightarrow \infty$ we get
\begin{equation}\label{eq:21}
\begin{split}
E(X_i^r) &\leq lim_{k\rightarrow\infty}\Big\{\sigma_i^{2(1-\frac{1}{2^k})}\Big\}\\
& \Rightarrow E(X_i^r)\leq \sigma_i^2
\end{split}
\end{equation}
Therefore, from Eq. (\ref{eq:19}) and Eq. (\ref{eq:20}) we get 
\begin{equation*}
F_i =\sum_{r=2}^\infty\frac{s^{r-2}E(X_i^r)}{r!\sigma_i^2} \leq \sum_{r=2}^\infty\frac{s^{r-2}\sigma_i^2}{r!\sigma_i^2}
\end{equation*}
Therefore,
\begin{equation*}
F_i \leq \frac{1}{s^2}\sum_{r=2}^\infty\frac{s^r}{r!}=\frac{1}{s^2}(e^s-1-s)
\end{equation*}
Applying this to Eq. (\ref{eq:20}) we get , 
\begin{equation} \label{eq:22}
E(e^{sX_i}) \leq e^{s^2\sigma_i^2\frac{1}{s^2}(e^s-1-s)}
\end{equation}
Using Chernoff Bound and Markov's inequality , we say that 
\begin{equation*}
\begin{split}
P[X\geq t] \leq e^{-st}E[e^{sX}]\\
\end{split}
\end{equation*}
where $X= X_1+X_2+...+X_n$ hence, 
\begin{equation*}
\begin{split}
P\bigg[\sum_{i=1}^nX_i> t\bigg] & \leq e^{-st}E\bigg[\prod_{i=1}^ne^{sX_i}\bigg]\\
& = e^{-st}\prod_{i=1}^nE[e^{sX_i}] \hspace{.8cm} \text{(as given independent rv's)}
\end{split}
\end{equation*}
Using Eq. (\ref{eq:22}) to this we get, 
\begin{equation*}
P\bigg[\sum_{i=1}^nX_i> t\bigg]\leq e^{-st}\prod_{i=1}^ne^{s^2\sigma_i^2\frac{1}{s^2}(e^s-1-s)}
\end{equation*}
As $\sigma^2 = \frac{1}{n}\sum_{i=1}^n\sigma_i^2 $, hence, 
\begin{equation}\label{eq:23}
\begin{split}
P\bigg[\sum_{i=1}^nX_i> t\bigg]&\leq e^{-st}e^{\sum_{i=1}^n{ \sigma_i^2(e^s-1-s)}}\\
&=e^{-st}e^{ n\sigma^2(e^s-1-s)} 
\end{split}
\end{equation}
now to obtain the closest bound we minimize R.H.S w.r.t $s$, therefore we get
\begin{equation*}
\begin{split}
\frac{de^{n\sigma^2(e^s-1-s)-st}}{ds}&=e^{n\sigma^2(e^s-1-s)-st}(n\sigma^2(e^s-1)-t)=0\\
&\Rightarrow e^s-1 = \frac{t}{n\sigma^2}
\end{split}
\end{equation*}
We get, 
\begin{equation*}
s=\log\bigg(1+\frac{t}{n\sigma^2}\bigg)
\end{equation*}
Using $s$ in Eq. (\ref{eq:23}), we have 
\begin{equation*}
\begin{split}
P\bigg[\sum_{i=1}^nX_i> t\bigg]&\leq e^{{-log(1+\frac{t}{n\sigma^2})t}+{ n\sigma^2(e^{log(1+\frac{t}{n\sigma^2})}-1-{log(1+\frac{t}{n\sigma^2})})}}\\
&=e^{{-log(1+\frac{t}{n\sigma^2})t}+{ n\sigma^2(\frac{t}{n\sigma^2}-{log(1+\frac{t}{n\sigma^2})})}}\\
&= e^{n\sigma^2(\frac{t}{n\sigma^2}-log(1+\frac{t}{n\sigma^2}) - \frac{t}{n\sigma^2}log(1+\frac{t}{n\sigma^2}))}
\end{split}
\end{equation*}
Let $h(u)=(1 + u)log(1 + u)-u $ for $u>0$, therefore we get
\begin{equation}\label{eq:24}
P\bigg[\sum_{i=1}^nX_i> t\bigg]\leq e^{-n\sigma^2h\Big(\frac{t}{n\sigma^2}\Big)} 
\end{equation}
\end{proof}
 \subsection{Bernstein's Inequality}
 \begin{tcolorbox}
 Under the same conditions defined in the Bennett's inequality, for any $\epsilon > 0,$ 
 \begin{center}
 $\mathbb{P}\bigg\{ \frac{1}{n}\sum_{i = 1}^n X_i > \epsilon \bigg\} \leq exp \bigg(- \frac{n\epsilon^2}{2(\sigma^2 + \epsilon/3)}\bigg)$
 \end{center}
 \end{tcolorbox}
\begin{proof}
 We can derive the Bernstein’s
inequality by further bounding the function h(x).
Let the function be, $G(x) =\frac{3}{2}\frac{x^2}{x+3}$. 
Now consider a function $\phi(x)=h(x)-G(x)$. $\phi''(x)= \frac{x^3+9x^2}{(x+1)(x+3)^3}$, For all $x\geq0$, $\phi\geq 0$ implies $\phi'(x)$ is increasing, i.e., for all $x\geq 0$, $\phi'(x)\geq 0$, and therefore $\phi(x)$ is increasing, hence $\phi(x)\geq0$ for all $x\geq0$, Hence we have
 
\begin{equation*}
h(x)\geq G(x)\ \ \forall x\geq 0
\end{equation*}
Therefore using Eq. (\ref{eq:24}) we get
\begin{equation*}
\begin{split}
&P\bigg[\sum_{i=1}^nX_i> t\bigg]\leq e^{-n\sigma^2G(\frac{t}{n\sigma^2})}\\
\Rightarrow &P\bigg[\sum_{i=1}^nX_i> t\bigg]\leq e^{(\frac{-3t^2}{2(t+3n\sigma^2)})}
\end{split}
\end{equation*}
Now let $t=n\epsilon$. Therefore,
\begin{equation}\label{eq:25}
P\bigg[\frac{1}{n}\sum_{i=1}^nX_i> \epsilon\bigg]\leq e^{-\frac{n\epsilon^2}{2\sigma^2+2\frac{\epsilon}{3}}} \hspace{1cm} 
\end{equation}
\end{proof}

\begin{example}
We have $n=2$ investments. Expected payoff of Investment 1 is $\mu_1 =\$50$ with standard deviation of $\sigma_1=\$25 $. Investment 2 has expected payoff $\mu_2 =\$70$ with standard deviation $\sigma_2=\$20$. Investment 1 has a floor on its payoff of $L_1=\$25$ and the upper bound of this payoff if $M_1 =\$65$. Meanwhile,  Investment 2 has it's floor payoff of $L_2 = \$60$ and ceiling payoff be $M_2=\$80$. For the portfolio to be worthwhile, we are told that the total payoff of both
investments must be at least \$130. We apply Bennett’s inequality, Bernstein’s
inequality and Hoeffding's inequality to this portfolio problem. If we calculate the probability bound using generic form of Bennett’s inequality
\begin{equation*}
P\bigg\{\frac{1}{n}(\sum_{i=1}^nX_i-\sum_{i=1}^nE[X_i])\geq t\bigg\}\leq \exp\bigg(\frac{-nv}{s^2}h\bigg(\frac{ts}{v}\bigg)\bigg)
\end{equation*}
where\\
$h(x)=(1+x)\ln(1+x)-x$\\
$s=\max_i(M_i-\mu_i)$\\
$v=\frac{1}{n}\sum_{i=1}^n\sigma_i^2$\\ \\
\text{The probability of complementary event specified in the inequality in turns out to be at least $0.9545$ for the values given in the example.}\\ \\
According to the generic form of Bernstein's inequality, 
\begin{equation*}
P\bigg\{(\sum_{i=1}^nX_i-\sum_{i=1}^nE[X_i])\geq t\bigg\} \leq \exp\bigg( \frac{-t^2}{2(n\sigma^2+(t/3))} \bigg)
\end{equation*}
where\\
$\sigma^2=\frac{1}{n}\sum_{i=1}^n\sigma_i^2$
\\Bernstein's gives the probability to be at least 0.9525.\\ \\
Applying Hoeffding's inequality to the same, we get
\begin{equation*}
P\bigg\{\sum_{i=1}^n(X_i-E[X_i])\geq t\bigg\}\leq \exp\bigg(\frac{-2t^2}{\sum_{i=1}^n(M_i-L_i)^2}\bigg)
\end{equation*}
where $M_i$ and $L_i$ are as specified in the example.\\
Hoeffding's gives the probability to be least 0.9048.\\ \\
Clearly Hoeffding's inequality gives the tightest bound in most of the cases.
\end{example}

 \subsection{Efron-Stein Inequality}
 \begin{tcolorbox}
 Let $\chi$ be some set and let $g : \chi^n \rightarrow \mathbb{R}$ be a  measurable
function of $n$ variables, $Z = g(X_1,.....,X_n)$ and its expected value is $\mathbb{E}(Z)$ where $X_1,....,X_n$ are arbitrary independent (not
necessarily identically distributed!) random variables taking values in $\chi$, Then
\begin{center}
$Var(Z) \leq \sum_{i=1}^n \mathbb{E}[(Z - \mathbb{E}_i(Z))^2]$
\end{center}
Where
\begin{equation*}
E_i(Z)=E[Z\mid X_1, X_2,...,X_{i-1}, X_{i+1},...,X_n]
\end{equation*}
\end{tcolorbox}
\begin{proof}
Let $V = Z-E(Z)$. Now if we define $V_i$ as
\begin{equation*}
V_i=E[Z|X_1,...,X_i]-E[Z|X_1,...,X_{i-1}] \ \ \forall i=2,...,n.
\end{equation*}
and for i=1,
\begin{equation*}
V_1=E[Z|X_1]-E[Z]
\end{equation*}
then 
\begin{equation*}
V=\sum_{i=1}^nV_i
\end{equation*}
and 
\begin{equation} \label{eq:26}
\begin{split}
Var(Z)&=E(V^2)\\
&=E\bigg(\bigg(\sum_{i=1}^nV_i\bigg)^2\bigg)\\
&=E\bigg(\sum_{i=1}^nV_i^2\bigg)+2E\bigg(\sum_{i>j}V_iV_j\bigg)
\end{split}
\end{equation}
now, $ E[XY]=E[E[XY|Y]]=E[YE[X|Y]] $ Therefore
\begin{equation}\label{eq:27}
E[V_iV_j]=E[V_jE[V_i|X_1,...,X_j]]
\end{equation}
Now we calculate 
\begin{equation*}
\begin{split}
E[V_i|X_1,...,X_j]&=E[(E[Z|X_1,...X_i]-E[Z|X_1,...,X_{i-1}])|X_1,...,X_j]\\
&=E[E[(Z|X_1,...X_i)|X_1,...,X_j]-E[(Z|X_1,...,X_{i-1})|X_1,...,X_j]]\\
\end{split}
\end{equation*}
Since $i>j$ and $i-1\geq j$ Then by Towering property
\begin{equation*}
E[V_i|X_1,...,X_j]= E[E[Z|X_1,...,X_j]-E[Z|X_1,...,X_j]]=0
\end{equation*}
Using this in Eq. (\ref{eq:27}) we get,
\begin{equation*}
E[V_iV_j]=0
\end{equation*}
Hence we have,
\begin{equation*}
Var(Z)=E\bigg(\sum_{i=1}^nV_i^2\bigg)=\sum_{i=1}^nE(V_i^2)
\end{equation*}
Bounding $E[V_i^2]$, 
\begin{equation*}
\begin{split}
V_i^2&=(E[Z|X_1,...,X_i]-E[Z|X_1,...,X_{i-1}])^2\\
&=(E[E[Z|X_1,...,X_n]-E[Z|X_1,...,X_{i-1},X_{i+1},...,X_n]|X_1,...,X_i])^2\\
&\leq E[(E[Z|X_1,...,X_n]-E[Z|X_1,...,X_{i-1},X_{i+1},...,X_n])^2|X_1,...,X_i]\\
& = E[(Z-E_i(Z))^2|X_1,...,X_i]
\end{split}
\end{equation*}
Summing over all $i$'s and taking expectation on both sides. As we know quadratic function is convex and hence we can apply Jensens inequality.
\begin{equation*}
Var(Z)\leq \sum_{i=1}^nE[(Z-E_i[Z])^2] 
\end{equation*}
\end{proof}

\begin{example}
\textbf {Kernel density estimation}\\
Let $X_1, . . . , X_n$ be i.i.d. real samples drawn according to some
density $\phi$. The kernel density estimate is
\begin{equation*}
\phi_n(x)=\frac{1}{nh}\sum_{i=1}^nK\bigg(\frac{x-X_i}{h}\bigg)
\end{equation*}
where $h > 0$ , and K is a nonnegative “kernel” $\int K = 1$. The $L_1$ error is
\begin{equation*}
Z=f(X_1,...,X_n)=\int|\phi(x)-\phi_n(x)|dx .
\end{equation*}
It is easy to see that
\begin{equation*}
|f(X_1,...,X_n)-f(X_1,...,X_i',...,X_n)|\leq \frac{1}{nh}\int\bigg| K\bigg(\frac{x-X_i}{h}\bigg)-K\bigg(\frac{x-X_i'}{h}\bigg)\bigg|\leq\frac{2}{n}
\end{equation*}
so we get 
\begin{equation*}
Var(Z)\leq \frac{2}{n}.
\end{equation*}
\end{example}
\subsection{McDiarmid's Inequality}
\begin{tcolorbox}
Let $X_1,....,X_m$ be independent random variables all taking values in the set $\chi$. Further, let $f : \chi^m \rightarrow \mathbb{R}$ be a function of $X_1,....,X_m$ that satisfies $\forall i, \forall X_1,...,X_m,X_i^{'} \in \chi$,
\begin{equation*}
|f(X_1,...,X_i,...,X_m) - f(X_1,...,X_i^{'},...,X_m)| \leq c_i
\end{equation*}
Then for all $\epsilon > 0$,
\begin{equation*}
Pr[f - \mathbb{E}[f] \geq \epsilon ] \leq exp\bigg( \frac{-2\epsilon^2}{\sum_{i=1}^m c_i^2}\bigg)
\end{equation*}
\end{tcolorbox}
\begin{proof}
Let $f'_i=f(X_1,...,X'_i,...,X_n)$\\
Using Chernoff bound we get,
\begin{equation*}
P[f-E[f]\geq \epsilon]\leq e^{-\epsilon s}e^{sE[f-E[f]]}
\end{equation*}
Now let,
\begin{equation*}
V_i=E[f|X_1,...,X_i]-E[f|X_1,...,X_{i-1}] \hspace{.5cm}\forall i=1,...,n
\end{equation*}
then 
\begin{equation*}
V=\sum_{i=1}^nV_i=f-E[f]
\end{equation*}
Therefore,
\begin{equation}\label{eq:28}
\begin{split}
P[f-E[f]\geq \epsilon]& \leq e^{-\epsilon s}E[e^{\sum_{i=1}^nsV_i}]\\
&=e^{-\epsilon s}\prod_{i=1}^nE[e^{sV_i}]
\end{split}
\end{equation}
Now let $ V_i$ be bounded by the interval $[L_i, U_i]$. We know that $|f-f'_i|\leq c_i$  ,
hence it follows that $|V_i|  \leq c_i$ and hence $|U_i-L_i| \leq c_i$. Using Hoeffding's lemma on $E[e^{sV_i} ]$
we get,
\begin{equation*}
E[e^{sV_i}]\leq e^{\frac{s^2(U_i-L_i)^2}{8}}\leq e^{\frac{s^2c_i^2}{8}}
\end{equation*}
Using this in Eq. (\ref{eq:28}) we get,
\begin{equation*}
\begin{split}
P[f-E[f]\geq\epsilon]&\leq e^{-\epsilon s}\prod_{i=1}^ne^{\frac{s^2c_i^2}{8}}\\
&=e^{-s\epsilon+s^2\sum_{i=1}^n\frac{c_i^2}{8}}
\end{split}
\end{equation*}
Now to make the bound tight we simply minimize it with respect to s.
Therefore,
\begin{equation*}
2s\sum_{i=1}^n\frac{c_i^2}{8}-\epsilon=0
\end{equation*}
\begin{equation*}
\Rightarrow s=\frac{4\epsilon}{\sum_{i=1}^nc_i^2}
\end{equation*}
Hence the bound is given by,
\begin{equation*}
\begin{split}
&P[f-E[f]\geq\epsilon]\leq e^{-\frac{4\epsilon}{\sum_{i=1}^nc_i^2}\epsilon+{(\frac{4\epsilon}{\sum_{i=1}^nc_i^2})}^2\sum_{i=1}^n\frac{c_i^2}{8}}\\
&\Rightarrow P[f-E[f]\geq\epsilon]\leq e^{-\frac{2\epsilon^2}{\sum_{i=1}^nc_i}} 
\end{split}
\end{equation*}
\end{proof}
\begin{example}
Let $X_1,...,X_n \in A$ be n-tuple i.i.d. random variables whose common distribution is $P$, i.e. $X_1,...,X_n \sim P $ and let $P_n(A)$ be the empirical distribution. The empirical distribution assigns the probability $1/n$ to each $X_i$
\begin{equation*}
P_n(A) =\frac{1}{n}\sum_{i=1}^nI(X_i \in A) 
\end{equation*} 
Define $\triangle_n \equiv f(X_1,...,X_n) = sup_A |P_n(A)
P(A)|$. Changing one observation changes $f$ by at most $\frac{1}{n}$. Hence,
\begin{equation*}
P\bigg(|\triangle_n-E(\triangle_n)|>\epsilon\bigg)\leq 2e^{-2n\epsilon^2}.
\end{equation*}\\
\end{example}
\begin{example}
\textbf{Kernel density function\\}
Similar to the example for Efron-Stein inequality, $X_1,X_2,...,X_n$ be i.i.d. random variable and $\phi_n(x)$ be the kernel density estimate. If $Z=f(X_1,...X_n)=\int|\phi(x)-\phi_n(x)|dx$
then, $|f(X_1,...,X_n)-f(X_1,...,X'_i,...,X_n)|\leq \frac{2}{n}$.
Thus, we can observe that $f (X_n)$ has the bounded differences property with $c_1 = ... = c_n = 2/n$. Applying McDiarmid's inequality on $f(x)$ we get, $P(|f(X_n)-E[f(X_n)]|\geq \epsilon)\leq 2e^{-n\epsilon^2/2}$.
\end{example}

\section{References}
\begin{itemize}
\item Mitzenmacher, Michael; Upfal, Eli (2005). Probability and Computing: Randomized Algorithms and Probabilistic Analysis. Cambridge University Press
\item Maxim Raginsky, Igal Sason, Concentration of Measure Inequalities in
Information Theory, Communications, and Coding.
\item Stephane Boucheron, Gabor Lugosi, Pascal Massart, Concentration inequalities
A nonasymptotic theory of independence
\item Anna Karlin, CSE525: Randomized Algorithms and Probabilistic Analysis, Lecture 18
\item Maxim Raginsky, Concentration inequalities
\end{itemize}
\end{document}